\documentclass[11pt,letterpaper,dvips]{article}
 \usepackage{amssymb,latexsym,amsmath,amsthm}
\newtheorem{thm*}{Theorem}

\newtheorem{thm}{Theorem}

\newtheorem{lemma}{Lemma}

\newtheorem{remark}{Remark}

\begin{document}

\def\d{ \partial_{x_j} } 
\def\Na{{\mathbb{N}}}

\def\Z{{\mathbb{Z}}}

\def\IR{{\mathbb{R}}}

\newcommand{\E}[0]{ \varepsilon}

\newcommand{\la}[0]{ \lambda}

\newcommand{\s}[0]{ \mathcal{S}}

\newcommand{\AO}[1]{\| #1 \| }

\newcommand{\BO}[2]{ \left( #1 , #2 \right) }

\newcommand{\CO}[2]{ \left\langle #1 , #2 \right\rangle} 

\newcommand{\R}[0]{ \IR\cup \{\infty \} } 

\newcommand{\co}[1]{ #1^{\prime}} 

\newcommand{\p}[0]{ p^{\prime}} 

\newcommand{\m}[1]{   \mathcal{ #1 }} 

\newcommand{ \W}[0]{ \mathcal{W}}

%  norm of H
\newcommand{ \A}[1]{ \left\| #1 \right\|_H }

% inner product H 
\newcommand{\B}[2]{ \left( #1 , #2 \right)_H }

% H^* , H pairing
\newcommand{\C}[2]{ \left\langle #1 , #2 \right\rangle_{  H^* , H } }

 \newcommand{\HON}[1]{ \| #1 \|_{ H^1} }

% Omega    \Om
\newcommand{ \Om }{ \Omega}

% \partial Omega      \pOm 
\newcommand{ \pOm}{\partial \Omega}

%   D(\Omega)   \D
\newcommand{\D}{ \mathcal{D} \left( \Omega \right)}

% D'( Omega)        \DP 
\newcommand{\DP}{ \mathcal{D}^{\prime} \left( \Omega \right)  }

% D' pairing
\newcommand{\DPP}[2]{   \left\langle #1 , #2 \right\rangle_{  \mathcal{D}^{\prime}, \mathcal{D} }}

% (H^1)^* , H^1    (( pairing ))      \PHH
\newcommand{\PHH}[2]{    \left\langle #1 , #2 \right\rangle_{    \left(H^1 \right)^*  ,  H^1   }    }

%  H^{-1} , H_0^1  (( pairing ))   \PHO
\newcommand{\PHO}[2]{  \left\langle #1 , #2 \right\rangle_{  H^{-1}  , H_0^1  }} 

 %  H^1(\Omega)     \HO
 \newcommand{\HO}{ H^1 \left( \Omega \right)}

%  H_0^1( \Omega)       \HOO
\newcommand{\HOO}{ H_0^1 \left( \Omega \right) }

% C_c^\infty(omega)
\newcommand{\CC}{C_c^\infty\left(\Omega \right) }

%H_0^1(Omega)  norm
\newcommand{\N}[1]{ \left\| #1\right\|_{ H_0^1  }  }

%H_0^1(Omega)   innerproduct 
\newcommand{\IN}[2]{ \left(#1,#2\right)_{  H_0^1} }

% H^1(\Omega) inner product
\newcommand{\INI}[2]{ \left( #1 ,#2 \right)_ { H^1}} 

% (H^1(\Omega))^*
\newcommand{\HH}{   H^1 \left( \Omega \right)^* } 

% ( H^{-1}(\Omega))
\newcommand{\HL}{ H^{-1} \left( \Omega \right) }

\newcommand{\HS}[1]{ \| #1 \|_{H^*}}

\newcommand{\HSI}[2]{ \left( #1 , #2 \right)_{ H^*}}

\newcommand{\WO}{ W_0^{1,p}} 
\newcommand{\w}[1]{ \| #1 \|_{W_0^{1,p}}}  

\newcommand{\ww}{(W_0^{1,p})^*}   

\newcommand{\Ov}{ \overline{\Omega}} 

\author{Craig Cowan}

\title{Regularity of the extremal solutions in a Gelfand system problem} 
\maketitle

\begin{abstract} We examine the elliptic system given by 
\begin{equation*} 
(P)_{\lambda,\gamma} \qquad \left\{
\begin{array}{rrl}
-\Delta u &=& \lambda e^v \qquad \Omega \\
-\Delta v &=& \gamma e^u \qquad \Omega \\
u &=& 0 \qquad \quad \pOm  \\
v&=& 0 \qquad \quad \pOm
\end{array}
\right.
\end{equation*} where $ \lambda, \gamma $ are positive parameters and where $ \Omega$ is a smooth bounded domain in $ \IR^N$.  Let $ \mathcal{U}$ denote the parameter region $ (\lambda, \gamma)$ of strictly positive parameters where $(P)_{\lambda,\gamma}$ has a smooth solution and let $ \Upsilon$ denote the boundary of $ \mathcal{U}$.  We show that the extremal solution $(u,v)$ associated with $ (\lambda, \gamma) \in \Upsilon$ is smooth provided that $ 3 \le N \le 9$ and 
\[  \frac{N-2}{8} < \frac{\gamma}{\lambda} < \frac{8}{N-2}.\]

\end{abstract}

\section{Introduction}

In this short note we are interested in  solutions of the  elliptic system given by 

\begin{equation*} 
(P)_{\lambda,\gamma} \qquad \left\{
\begin{array}{rrl}
-\Delta u &=& \lambda e^v \qquad \Omega \\
-\Delta v &=& \gamma e^u \qquad \Omega \\
u &=& 0 \qquad \quad \pOm  \\
v&=& 0 \qquad \quad \pOm
\end{array}
\right.
\end{equation*} where $ \lambda, \gamma $ are positive parameters and where $ \Omega$ is a smooth bounded domain in $ \IR^N$.          In particular we are interested in the regularity of the extremal solutions associated with $(P)_{\lambda,\gamma}$, which we define more precisely later.   Along the diagonal $ \lambda = \gamma$ the problem $ (P)_{\lambda,\gamma}$ reduces to the scalar analog of $ (P)_{\lambda,\gamma}$, see below. Provided one stays sufficiently close to the diagonal we show that some basic maximum principle arguments coupled with  a standard energy estimate approach (the familiar approach in the scalar case) shows the regularity of the extremal solutions in the expected dimensions.    

    We now recall  the well studied scalar version (with general nonlinearity $f$)  of $(P)_{\lambda,\gamma}$  given by 
  \begin{equation*} 
(P)_\lambda \qquad \left\{
\begin{array}{rrl}
-\Delta u &=& \lambda f(u) \qquad \Omega \\
u &=& 0 \qquad \quad \pOm  \\
\end{array}
\right.
\end{equation*}   where $ \lambda$ is a positive parameter and where $ \Omega$ is a bounded domain in $ \IR^N$.  See, for instance,    \cite{BCMR}, \cite{BV}, \cite{CR}, \cite{M} and \cite{MP}.
Here generally one assumes that $ f$ is a smooth, increasing, convex nonlinearity with $ f(0)=1$ and $ f$ superlinear at $ \infty$, ie. $ \lim_{u \rightarrow \infty} \frac{f(u)}{u} = \infty$. 
 It is known that there is an non degenerate finite interval $ \mathcal{U} = (0,\lambda^*) $ such that for all $ 0 < \lambda < \lambda^*$ there exists a smooth, \textbf{minimal solution} $ u_\lambda$ of $(P)_\lambda$.  By minimal we mean that any other solution $ v$ of $ (P)_\lambda$ satisfies $ v \ge u_\lambda$ a.e. in $ \Omega$.   In addition one can show that for each $ x \in \Omega$ the map $ \lambda \mapsto u_\lambda(x)$ is increasing on $(0,\lambda^*)$.  This allows one to define the \textbf{extremal solution}  
\[ u^*(x):= \lim_{\lambda \nearrow \lambda^*} u_\lambda(x),\]  and it can be shown that $ u^*$ is the unique weak solution of $ (P)_{\lambda^*}$.  Also it is known that for $ \lambda > \lambda^*$ there are no weak solutions.   One can also show that for each $ 0 < \lambda < \lambda^*$ the minimal solution $ u_\lambda$ is semi-stable in the sense that the principle eigenvalue of the linear operator 
\[ L_{\lambda,u_\lambda}:= -\Delta - \lambda f'(u_\lambda),\] over $H_0^1(\Omega)$ is nonnegative.    Using the variational structure this implies that 
\[ \int_\Omega \lambda f'(u_\lambda) \psi^2 dx \le \int_\Omega | \nabla \psi|^2 dx, \qquad \forall \psi \in H_0^1(\Omega).\]    One can now ask the question whether $ u^*$ is  a classical solution of $(P)_{\lambda^*}$?  Elliptic regularity shows this is equivalent to the boundedness of $ u^*$.     In the case where $ f(u) = e^u$ one can show  that $ u^*$ is bounded provided $ N \le 9$.  Moreover this is optimal after one considers the fact that $ u^*(x)= -2 \log(|x|)$ provided $ \Omega$ is the unit ball in $ \IR^N$ where $ N \ge 10$.      For more results concerning the regularity of the extremal solution $ u^*$ the reader should see  \cite{Nedev}, \cite{CC}, \cite{Cabre} and \cite{Ye}.      We mention that vital to all the results concerning the regularity of $u^*$ is to use the semi-stability of the minimal solutions $ u_\lambda$ to obtain a priori estimates and then to pass to the limit.

We now return to the system $(P)_{\lambda,\gamma}$ and  we follow the work of M. Montenegro  \cite{Mont}, where all of the following results are taken from.     We also mention that he obtains many more results and also that he studies a much more general system then $(P)_{\lambda,\gamma}$.  We let  $ \mathcal{Q}=\{ (\lambda,\gamma): \lambda, \gamma >0 \}$ and we define 
\[ \mathcal{U}:= \left\{ (\lambda,\gamma) \in \mathcal{Q}: \mbox{ there exists a smooth solution $(u,v)$ of $(P)_{\lambda,\gamma}$} \right\}.\]

We set $ \Upsilon:= \partial \mathcal{U} \cap \mathcal{Q}$.   The curve $ \Upsilon$ is well defined and separates $ \mathcal{Q}$ into two connected components $ \mathcal{Q}$ and $ \mathcal{V}$.   We omit the various properties of $ \Upsilon$ but the interested reader should consult \cite{Mont}.   One point we mention is that if for $ x,y \in \IR^2$ we say $ x \le y$ provided $ x_i \le y_i$ for $ i=1,2$ then it is easily seen, using the method of sub/supersolutions, that if $ (0,0) < (\lambda_0,\gamma_0) \le (\lambda,\gamma) \in \mathcal{U}$ then $ (\lambda_0,\gamma_0) \in \mathcal{U}$.    Now it can be shown that for each $ (\lambda,\gamma) \in \mathcal{U}$ there exists a smooth minimal solution $(u_{\lambda,\gamma},v_{\lambda,\gamma})$ of $ (P)_{\lambda,\gamma}$ and  if $(0,0)< (\lambda_1,\gamma_1) \le (\lambda_2,\gamma_2) \in \mathcal{U}$  then 
\[ ( u_{\lambda_1,\gamma_1}, v_{\lambda_1,\gamma_1}) \le (u_{\lambda_2,\gamma_2}, v_{\lambda_2,\gamma_2}).\]   Now for $ (\lambda^*,\gamma^*) \in \Upsilon$  there is some $ 0 < \sigma < \infty$ such that $ \gamma^* = \sigma \lambda^*$ and we can define the extremal solution $(u^*,v^*)$ at $ (\lambda^*,\gamma^*)$ by passing to the limit along the ray given by $ \gamma = \sigma \lambda$ for $ 0 < \lambda < \lambda^*$.   Moreover it can be shown that $ (u^*,v^*)$ is indeed a weak solution of $ (P)_{\lambda^*,\gamma^*}$.            We now come to the issue of stability.     

\begin{thm} \label{stable} (\cite{Mont})   Let $ (\lambda,\gamma) \in \mathcal{U}$ and let $ (u,v)$ denote the minimal solution of $(P)_{\lambda,\gamma}$.  Then $(u,v)$ is semi-stable in the sense that there is some smooth $ 0 < \phi,\psi \in H_0^1(\Omega)$ and $  0 \le K $ such that 
\[ -\Delta \phi = \lambda e^v \psi + K \phi, \qquad -\Delta \psi = \gamma e^u \phi + K \psi, \qquad \Omega.\] 

\end{thm}   Now one should note that $ K < \lambda_1(\Omega)$.  To see this one multiplies either of these equations by the first positive eigenfunction of $ -\Delta$ and integrates by parts. 

\section{Main Results}

Our main result is given by the following theorem.

\begin{thm} \label{main} Let $ 3 \le N \le 9$ and suppose that $ (\lambda,\gamma) \in \Upsilon$ with 
\[ \frac{N-2}{8} < \frac{\gamma}{\lambda} < \frac{8}{N-2}.\]    Then the associated extremal solution $ (u^*,v^*)$ is smooth. 

\end{thm}

One should note that along the diagonal the problem reduces to the scalar problem.  Also by symmetry it is enough to prove the result for $ 0 < \gamma \le \lambda$.   We prove the above Theorem in a series of lemma's.

\begin{lemma} Suppose that $ (u,v)$ is a smooth solution of $ (P)_{\lambda,\gamma}$ where $ 0 < \gamma \le \lambda$.  Then 
\[ \frac{\gamma}{\lambda} u \le v \le u \qquad \mbox{ a.e. in $ \Omega$.}\]

\end{lemma}  

\begin{proof} Taking the difference of the equations in $ (P)_{\lambda,\gamma}$ we have  $ -\Delta (u-v) = \lambda e^v - \gamma e^u$ in $ \Omega$ and multiplying this by $ (u-v)_-$ and integrating by parts one arrives at 
\[ - \int_\Omega | \nabla (u-v)_-|^2 dx  = \int_\Omega ( \lambda e^v - \gamma e^u) (u-v)_- dx,\]  and now note that the right hand side is nonnegative where as the left hand side is nonpositive.  Hence we see that $ (u-v)_-=0$ a.e. in $ \Omega$ and so $ u \ge v$ a.e. in $ \Omega$.     Now note that 
\[ -\Delta (v-\frac{\gamma}{\lambda}u) = \gamma (e^u - e^v) \ge 0 \qquad \Omega,\] since $ u \ge v$ in $ \Omega$ and so $ v \ge \frac{\gamma}{\lambda} u$ in $ \Omega$.

\end{proof}

\begin{lemma} Suppose that $(\lambda,\gamma) \in \mathcal{U}$ with $ 0 < \gamma \le \lambda$ and we let $(u,v)$ denote the minimal solution of $(P)_{\lambda,\gamma}$.  Let $K,\phi,\psi$ be as in Theorem \ref{stable}.  Then 
\[ \frac{\psi}{\phi} \ge \frac{\gamma}{\lambda} \qquad \mbox{ in $\Omega$.}\]

\end{lemma}

\begin{proof}

 First note that 
\begin{eqnarray*}
-\Delta ( \psi - \phi) &=& \gamma e^u \phi - \lambda e^v \psi + K( \psi - \phi) \\
& \ge & \gamma e^v ( \phi - \psi) +( \gamma - \lambda) e^v \psi + K( \psi-\phi)
\end{eqnarray*} where we have used the fact that $  u \ge v$ in $ \Omega$. Rearranging this we have 
\begin{equation} \label{eq1}
-\Delta ( \psi - \phi) - K( \psi-\phi) + \gamma e^v (\psi - \phi) \ge ( \gamma- \lambda) e^v \psi \qquad \Omega.
\end{equation} We now define $L:=-\Delta -K$.  One now notes that 
\begin{eqnarray*}
L( \psi - \phi + \frac{ \lambda-\gamma}{\lambda} \phi) + \gamma e^v ( \psi - \phi + \frac{ \lambda-\gamma}{\lambda} \phi) & \ge & L( \psi - \phi + \frac{ \lambda-\gamma}{\lambda} \phi) + \gamma e^v ( \psi - \phi) \\
&=& L( \psi - \phi) + \gamma e^v(\psi - \phi) \\ 
&=& L( \psi-\phi) + \gamma e^v(\psi - \phi) + \frac{ \lambda- \gamma}{\lambda} L(\phi) \\
& \ge & (\gamma - \lambda) e^v \psi + \frac{ \lambda- \gamma}{\lambda} L(\phi) \\
&=& (\gamma - \lambda) e^v \psi + \frac{ \lambda- \gamma}{\lambda} \lambda e^v \psi \\
&=&0 
\end{eqnarray*} in $ \Omega$.  Now since $ L + \gamma e^v$ satisfies the maximum principle we see that $ \psi - \phi + \frac{\lambda- \gamma}{\lambda} \phi \ge 0$ in $ \Omega$, which after re-arranging, gives the desired result.

\end{proof}

Theorem \ref{main} will easily follow from the following lemma. 

\begin{lemma} Suppose $ 3 \le N \le 9$ and that $(u_m,v_m)$ denotes a sequence of smooth minimal solutions to $(P)_{\lambda_m, \sigma \lambda_m}$ where $ \frac{N-2}{8} < \sigma \le 1$.  Then $ (u_m,v_m)$ is bounded in $L^\infty(\Omega) \times L^\infty(\Omega)$.  

\end{lemma}

\begin{proof}  Fix $  \frac{N-2}{8} < \sigma \le 1$ and for notational simplicity we drop the subscript $m$ from $ u_m,v_m,\phi_m,\psi_m$ and $K_m$.  From the previous lemma we have $ \frac{ \psi}{\phi} \ge \sigma$.  Now for any smooth positive function $ E$ one has 
\[ \int_\Omega \frac{-\Delta E}{E} \beta^2 dx \le \int_\Omega | \nabla \beta|^2 dx, \qquad \forall \beta \in H_0^1(\Omega),\] see, for instance, \cite{craig}.

  Now note that 
\[ \frac{-\Delta \phi}{\phi} = \lambda e^v \frac{\psi}{\phi} + K \ge \sigma \lambda e^v \qquad \Omega.\] Taking $E=\phi$ and $ \beta = e^{tu}-1$, where $ t$ is chosen such that $ \frac{N-2}{4}<t< 2 \sigma$, gives 
\begin{equation} \label{z}
 \sigma \lambda \int_\Omega e^v \left( e^{tu} -1 \right)^2 dx \le t^2 \int_\Omega e^{2tu} | \nabla u|^2 dx.
 \end{equation}
  Now multiplying $ -\Delta u = \lambda e^v$ by $ e^{2tu}-1$ and integrating by parts gives 
\begin{equation} \label{zz}
 2 t \int_\Omega e^{2tu} | \nabla u|^2 dx = \lambda \int_\Omega e^v (e^{2tu}-1) dx.
 \end{equation}  Now equating (\ref{z}) and (\ref{zz}) gives, after some simplification, 
 \[ \left( \frac{\sigma}{t} - \frac{1}{2} \right) \int_\Omega e^v e^{2tu} dx \le \frac{2 \sigma}{t} \int_\Omega e^{tu} e^v dx.\]  Now note that since $ t < 2 \sigma$ the coefficient on the left is positive. Now applying Holder's inequality on the right and squaring gives 
 \[ \left( \frac{ \sigma}{t} - \frac{1}{2} \right)^2 \int_\Omega e^{2tu}e^v dx \le \frac{4 \sigma^2}{t^2} \int_\Omega e^v dx,\]  and now since $ u \ge v$ in $ \Omega$ we see that this gives us an $L^{2t+1}(\Omega)$ bound for $ e^v$.  We now return to the sequence notation.   So we have that $ e^{v_m}$ is bounded in $L^{2t+1}(\Omega)$ but note that $ 2t+1 > \frac{N}{2}$ and also note that $ \lambda_m $ is bounded.  Now since $ -\Delta u_m = \lambda_m e^{v_m}$ in $ \Omega$ with $ u_m =0$ on $ \pOm$, and since $ \lambda_m$ is bounded one sees, using elliptic regularity, that $ u_m$ is bounded in $ L^\infty(\Omega)$.   From this, and since $ \sigma \lambda_m$ is bounded, we easily infer that $ v_m$ is bounded in $L^\infty(\Omega)$. 
\end{proof}

\begin{remark} A natural system to examine is 
\begin{equation*} 
 \qquad \left\{
\begin{array}{rrl}
-\Delta u &=& \lambda (v+1)^p \qquad   \Omega \\
-\Delta v &=& \gamma  (u+1)^q \qquad   \Omega \\
u &=& 0 \qquad \qquad  \quad \pOm  \\
v&=& 0 \qquad \qquad  \quad \pOm
\end{array}
\right.
\end{equation*} where $ 1 < p,q$.  In the special case where $ p=q$ our methods easily gives similiar type results concerning the regularity of the extremal solutions.

\end{remark}

\end{document}